\documentclass[10pt]{article}

\usepackage{amsmath,amssymb,amsthm}

\usepackage{graphicx}
\usepackage{tikz}

\newtheorem*{theorem*}{Theorem}

\newcommand{\si}{\sim}
\newcommand{\bx}{\boxed}
\renewcommand{\do}{\dots}
\newcommand{\RR}{\bf R}
\newcommand{\CC}{\bf C}
\newcommand{\ZZ}{\bf Z}
\newcommand{\NN}{\bf N}
\newcommand{\rk}{\rm rk}
\newcommand{\rar}{\rightarrow}
\newcommand{\lra}{\leftrightarrow}
\renewcommand{\a}{\alpha}

\newcommand{\g}{\gamma}
\newcommand{\cb}{{\cal B}}
\newcommand{\sub}{\subset}
\newcommand{\n}{\notin}
\newcommand{\T}{\times}
\newcommand{\e}{\epsilon}

\newcommand{\co}{{\cal O}}
\newcommand{\m}{\mapsto}

\newcommand{\emp}{\emptyset}
\newcommand{\boc}{\bf c}

\newcommand{\codim}{\rm codim}

\title {On conjugacy classes in the Lie group $E_8$}
\author{George Lusztig\thanks{Supported in part by National Science Foundation grant
  DMS-6927854.}}
\date{}

\begin{document}
\maketitle

\vspace{-2ex}

\centerline{Based on a talk given at the AMS-RMS joint meeting in Alba Iulia, 6/28/2013.}
\vspace{4ex}

The classification of simple Lie algebras (or Lie groups)/$\CC$ (Killing 1879, a glory of 19-th century 
mathematics) leads to $4$ $\infty$ series and $5$ exceptional groups of which the
largest one, $E_8$, has dimension $248$. 

The quantity $\frac{\dim(G)}{\rk(G)^2}$ defined for any simple Lie group $G$ is bounded. It reaches its
maximum $\frac{248}{8^2}\cong4$) for $G$ of type $E_8$. This shows that the group of type $E_8$
is the most noncommutative of all simple Lie groups.

The Lie group $E_8$ can be obtained from the graph $E_8$ (to be described below)
by a method of Chevalley (1955), simplified using theory of "canonical bases" (1990).

The graph $E_8$ has $8$ vertices and $7$ edges. It can be obtained by removing $5$ edges
from the $1$-skeleton of a cube in $\RR^3$ (wich has $8$ vertices and $12$ edges).

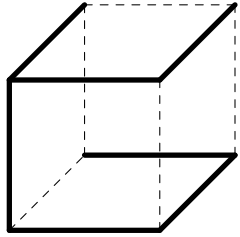
\begin{center}
\begin{tikzpicture}
\draw [thick/.style={line width=2pt}, %this number controls edge thickness
       line cap=round]                %rounds off the ends of the edges
 (0,0) edge [thick]  (2,0)
 (2,0) edge [dashed] (2,2)
 (2,2) edge [thick]  (0,2)
 (0,2) edge [thick]  (0,0)

 (1,1) edge [thick]  (3,1)
 (3,1) edge [dashed] (3,3)
 (3,3) edge [dashed] (1,3)
 (1,3) edge [dashed] (1,1)

 (0,0) edge [dashed] (1,1)
 (2,0) edge [thick]  (3,1)
 (2,2) edge [thick]  (3,3)
 (0,2) edge [thick]  (1,3) 
; %an important semicolon
\end{tikzpicture}
\end{center}

Let $I$ be the set of vertices of the graph $E_8$. Let $V$ be the free $\ZZ$-module with basis 
$\{\a_i;i\in I\}$ and bilinear form:

$(\a_i,\a_i)=2$; $(\a_i,\a_j)=-\sharp(\text{edges joining }i,j)$ if $i\ne j$.

This bilinear form was defined by Korkine-Zolotarev (1873) before 
Killing and nonconstructively earlier by Smith (1867). The set $R=\{\a\in V,(\a,\a)=2\}$ (roots)
has exactly $240$ elements. Its subset $R^+=R\cap\sum_{i\in I}\NN\a_i$ has exactly $120$ elements. 

Let $M$ be the $\CC$-vector with basis $X_\a (\a\in R)$, $t_i (i\in I)$. 
For $i\in I$, $\e=1,-1$, we define linear maps $E_{i,\e}:M\rar M$ by 
 
$E_{i,\e}X_\a=X_{\a+\e\a_i}$ if $\a\in R,\a+\e\a_i\in R$,

$E_{i,\e}X_\a=0$ if $\a\in R,\a+\e\a_i\n R\cup 0$,

$E_{i,\e}X_{-\e\a_i}=t_i$, $E_{i,\e}t_j=|(\a_i,\a_j)|X_{\e\a_i}$.
\newline
Let $M^+=\sum_{\a\in R^+}\CC X_\a\sub M$. We have $\dim M=248$, $\dim M^+=120$. 

Let $G=E_8(\CC)$ be the subgroup of $GL(M)$ generated by the automorphisms 
$\exp(\lambda E_{i,\e})$, $i\in I$, $\e=1,-1$, $\lambda\in\CC^*$. (The Lie group $E_8$.) 
Replacing $\CC$ by a finite field $F_q$ leads to a finite group $E_8(F_q)$ which, by Chevalley 
(1955), is simple of order $q^{248}+$ lower powers of $q$. (The denominators
in the exponential do not cause troubles.)

The Weyl group is $W=\{g\in Aut(V);gR=R\}$. It can be viewed as a limit of $E_8(F_q)$ as $q\to1$
hence we can denoted it by $E_8(F_1)$. It is a finite group:

$\sharp(E_8(F_1))=\lim_{q\to1}\frac{\sharp(E_8(F_q))}{(q-1)^8}=4!6!8!$.

For a group $H$, the {\it conjugacy classes} of $H$ are the orbits of 
the action $x:g\m xgx^{-1}$ of $H$ on itself. They form a set $cl(H)$.

The main guiding principle of this talk: 

{\it The conjugacy classes of $G=E_8(\CC)$ should be organized 
according to the conjugacy classes of $W=E_8(F_1)$.}

For $w\in W$ let $l(w)=\sharp(R^+\cap w(R-R^+))\in\NN$ (length).

For $C\in cl(W)$ let $C_{min}=\{w\in C;l:C\rar\NN\text{ reaches minimum at }w\}$.

$C\in cl(W)$ is {\it elliptic} if $\{x\in V;wx=x\}=0$ for some/any $w\in C$.

$cl(W)$ has been described by Carter (1972); we have
$\sharp cl(W)=112$, $\sharp\{C\in cl(W);C\text{ elliptic}\}=30$.

The $G$-orbit of $M^+$ in the Grassmannian of $120$-dimensional 
subspaces of $M$ is a closed smooth subvariety $\cb$ of dimension $120$, the {\it flag manifold}.
The diagonal $G$-action on $\cb\T\cb$ has orbits in canonical bijection 

$\co_w\lra w$ 

with $W$. (Bruhat 1954, Harish Chandra 1956).

For $w\in W$ let $G_w=\{g\in G;(B,g(B))\in\co_w\text{ for some }B\in\cb\}$.

For $C\in cl(W)$ let $G_C=G_w$ where $w\in C_{min}$; one shows (using results of
Geck-Pfeiffer 1993, simplified by He-Nie 2012) that $G_C$ is independent of the choice of $w$ in 
$C_{min}$; also, $G_C\ne\emp$, $G_C$ is a union of conjugacy classes. Let 

$\delta_C=\min_{\g\in cl(G);\g\sub G_C}\dim\g$,

$\bx{G_C}=\cup_{\g\in cl(G);\g\sub G_C,\dim\g=\delta_C}\g$.

Then $\bx{G_C}$ is $\ne\emp$, a union of conjugacy classes of fixed dimension, $\delta_C$.

Let $G^{un}=\{g\in G;\text{ all eigenvalues of }g:M\rar M\text{ are }1\}$.
\vspace{1ex}

%\noindent{\bf Theorem} 
%\vspace{1ex}

\begin{theorem*}

\begin{list}{}{}
\item (a) $\cup_{C\in cl(W)}\bx{G_C}=G$.  
\item (b) For any $C,C'\in cl(W)$, $\bx{G_C},\bx{G_{C'}}$ are either equal or disjoint. 
\item (c) For any $C\in cl(W)$, $\bx{G_C}\cap G^{un}$ is either empty or a single conjugacy class.  
\item (d) If $C\in cl(W)$ is elliptic then $\delta_C=248-l(w)$ for any $w\in C_{min}$ and $\bx{G_C}\cap G^{un}$ 
is a single conjugacy class.   

\end{list}

%\endproclaim
%\vspace{1ex}

\end{theorem*}

See arxiv:1305.7168. 
\vspace{1ex}

The proof uses representation theory of $E_8(F_q)$ (a part of its 
character table known since 1980's) and  computer calculation. (Help with programming in GAP was provided 
by Gongqin Li.)

The subsets $\bx{G_C}$ partition $G$ into the {\it strata} of $G$.
If $\g,\g'\in cl(G)$ we say: $\g\si\g'$ if $\g,\g'$ are contained in the 
same stratum. If $C,C'\in cl(W)$ we say: $C\si C'$ if $\bx{G_C}=\bx{G_{C'}}$.
These are equivalence relations on $cl(G)=cl(E_8(\CC))$,
$cl(W)=cl(E_8(F_1))$; by the theorem we have canonically 

$cl(E_8(\CC))/\si\lra cl(E_8(F_1))/\si$.

{\it Examples.} If $C=\{1\}$ then $\bx{G_C}=\{1\}$.

If $C\in cl(W)$ contains all $w$ in $W$ of length $1$ then $\bx{G_C}$ is a single 
conjugacy class (it has dimension $58$.)

Let $C_{cox}=\{w\in W; w\text{ has order }30\}$; it is a single conjugacy 
class in $W$, the Coxeter class. If $C=C_{cox}$ then $G_C=\bx{G_C}$ is the 
union of all conjugacy classes of dimension $240$ (Steinberg 1965).

For $n=2,3,5$ let $\boc_n$ be the unique conjugacy class in $G$ such 
that $\codim\boc_n=240/n$ and any element of $\boc_n$ has order $n$.
One can show: $\boc_n\sub\bx{G_{C_{cox}^{30/n}}}$.

Let $cl_u(E_8(\CC))=\{\g\in cl(E_8(\CC);\g\sub G^{un}\}$. Let $cl_u(E_8(\bar F_p))$ be the
analogous set with $\CC$ replaced by $\bar F_p$, $p$ a prime number. We have
$\sharp(cl_u(E_8(\CC)))=70$, (Dynkin, Kostant), $\sharp(cl_u(E_8(\bar F_p)))=70+n$ 
where $n=4,1,0,0,\do$ for $p=2,3,5,7,\do$ (Mizuno 1980).

We have a natural imbedding $j_p:cl_u(E_8(\CC))\rar cl_u(E_8(\bar F_p))$ and
$cl(E_8(\CC))/\si\lra\cup_{p\text{ prime}}cl_u(E_8(\bar F_p))$ (union taken using $j_p$).
Hence the number of strata is $75$.

The definitions above extend to any simple Lie group $G$.

If $G=S0_5(\CC)$, then $cl(W)=\{C_4,C_4^2,C',C'',\{1\}\}$ where $C_4$ consists of the elements of order $4$;

$\bx{G_{C_4}}$ is the union of classes of dimension $8$, 

$\bx{G_{C_4^2}}$ is the union of classes of dimension $6$, 

$\bx{G_{C'}},\bx{G_{C''}}$ are the two conjugacy classes of dimension $4$, 

$\bx{G_{\{1\}}}=\{1\}$. 

\end{document}